\documentclass[12pt]{amsart}
\usepackage{amsfonts}
\usepackage{amssymb}
\usepackage{graphicx}

\def\D{{\mathbb D}}  
\def\C{{\mathbb C}}

\def\({\left(}       \def\){\right)}

\newtheorem{lem}{\sc Lemma}
\newtheorem{thm}{\sc Theorem}

\newtheorem{other}{\sc Theorem}              
\newenvironment{pf}{\noindent{\textit{Proof. }}}{$\Box$ }

\begin{document}
\title[Criteria for bounded valence of harmonic mappings]
{Criteria for bounded valence of harmonic mappings}

\author[J.-M. Huusko]{Juha-Matti Huusko}
\address{Department of Physics and Mathematics, University of Eastern Finland, P.O. Box 111, FI-80101 Joensuu, Finland.} \email{juha-matti.huusko@uef.fi}
\author[M. J. Mart\'{\i}n]{Mar\'{\i}a J. Mart\'{\i}n}
\address{Department of Physics and Mathematics, University of Eastern Finland, P.O. Box 111, FI-80101 Joensuu, Finland.} \email{maria.martin@uef.fi}

\date{\today}
\thanks{This research is supported in part by Academy of Finland grant $\#268009$. The first author is also supported by the Faculty of Science and Forestry of the University of Eastern Finland research project $\#930349$. The second author thankfully acknowledges partial support from grant Fondecyt $\#1150284$, Chile. Also, from Spanish MINECO/FEDER research project MTM2015-65792-P and by the Thematic Research Network MTM2015-69323-REDT, MINECO, Spain.}
\maketitle

\begin{abstract}
In 1984, Gehring and Pommerenke proved that if the Schwarzian derivative $S(f)$ of a locally univalent analytic function $f$ in the unit disk satisfies that $\limsup_{|z|\to 1} |S(f)(z)| (1-|z|^2)^2 < 2$, then there exists a positive integer $N$ such that $f$ takes every value at most $N$ times. Recently, Becker and Pommerenke have shown that the same result holds in those cases when the function $f$ satisfies that $\limsup_{|z|\to 1} |f''(z)/f'(z)|\, (1-|z|^2)< 1$.
\par
In this paper, we generalize these two criteria for bounded valence of analytic functions to the cases when $f$ is merely harmonic.
\end{abstract}

\maketitle
\date{\today}
\section*{Introduction}
\label{intro}
Let $\D$ be the unit disk in the complex plane $\mathbb C$. It is well known that if a locally univalent function $f$ in $\D$ satisfies
\[
\|P(f)\|=\sup_{z\in\D} \left|\frac{f''(z)}{f'(z)}\right|(1-|z|^2)\leq 1\,,
\]
then $f$ is globally univalent in $\D$. This criterion of univalence in due to Becker \cite{Becker}. Becker and Pommerenke showed that the constant $1$ is sharp \cite{Becker-Pom-1}.
\par\smallskip
The quotient $P(f)=f''/f'$ is the \emph{pre-Schwarzian derivative} of $f$. The quantity $\|P(f)\|$ defined above is said to be the \emph{pre-Schwarzian norm} of $f$.
\par\smallskip
Nehari \cite{Nehari} proved  that if a locally univalent analytic function $f$ in $\D$ satisfies
\begin{equation}\label{eq-Schwarzian}
\|S(f)\|=\sup_{z\in\D} |S(f)(z)|\, (1-|z|^2)^2\leq 2\,,
\end{equation}
then $f$ is globally univalent in $\D$. Here, $S(f)$ denotes the \emph{Schwarzian derivative} of $f$ defined by
\begin{equation}\label{eq-Schw-defn}
S(f)=P(f)'-\frac 12 (P(f))^2=\left(\frac{f''}{f'}\right)'-\frac 12 \left(\frac{f''}{f'}\right)^2\,.
\end{equation}
\par
The \emph{Schwarzian norm} $\|S(f)\|$ of $f$ equals the supremum in (\ref{eq-Schwarzian}).
\par\smallskip
The \emph{valence} of an analytic mapping $f$ in $\mathbb D$ is defined by $\sup_{w\in\C} n(f,w)$, where $n(f,w)$ is the number of points $z\in\D$ (counting multiplicities) for which $f(z)=w$. The function $f$ is said to have \emph{bounded valence} if there exists a positive integer $N$ such that $\sup_{w\in\C} n(f,w)\leq N$. That is, if there is a positive integer $N$ such that $f$ takes every value at most $N$ times in $\D$.
\par\smallskip
A criterion for the bounded valence of analytic functions in terms of the Schwarzian derivative has been known for some time. Binyamin Schwarz \cite{Schwarz}, using techniques from the theory of differential equations, proved that if a locally univalent analytic function $f$ in $\D$ satisfies
\[
|S(f)(z)|\, (1-|z|^2)^2\leq 2
\]
for all $z$ in an annulus $0\leq r_0< |z|<1$, then $f$ has bounded valence. The authors in \cite{GP} show that the slightly stronger condition stated in Theorem~\ref{thm-GP} below suffices to ensure not only that the locally univalent analytic function~$f$ in the unit disk has a spherically continuous extension to $\overline\D$ but also the criterion for bounded valence of analytic functions that we now enunciate.
\begin{other}\label{thm-GP}
Let $f$ be a locally univalent analytic function in the unit disk. If
\[
\limsup_{|z|\to 1} \left|S(f)(z)\right|(1-|z|^2)^2< 2\,,
\]
then $f$ has bounded valence.
\end{other}
\par
Only recently the corresponding bounded valence criterion to that stated in Theorem~\ref{thm-GP}, now in terms of the pre-Schwarzian derivative, has been obtained \cite[Thm.~3.4]{Becker-Pom-2}.
\begin{other}\label{thm-BP}
Let $f$ be a locally univalent analytic function in the unit disk. If
\[
\limsup_{|z|\to 1} \left|\frac{f''(z)}{f'(z)}\right|(1-|z|^2)< 1\,,
\]
then there exists a positive integer $N$ such that $f$ takes every value at most $N$ times in $\D$.
\end{other}
\par\smallskip
The main aim of this paper is to generalize these criteria stated in Theorems~\ref{thm-GP} and \ref{thm-BP} for bounded valence of locally univalent analytic functions in the unit disk to the cases when the function $f$ is merely \emph{harmonic}.
\par
Perhaps, it is appropriate to stress that we have not been able to find any paper containing bounded valence criteria for harmonic functions in $\mathbb D$. The article \cite{CDO}, which gathers bounded valence criteria for Weierstrass-Enneper \emph{lifts} of planar harmonic mappings to their associated minimal surfaces, deserves to be mentioned at this point.

\section{Background}
\par\smallskip
\subsection{Harmonic mappings}
A complex-valued harmonic function $f$ in the unit disk $\D$ can be written as $f=h+\overline g$, where both $h$ and $g$ are analytic in $\D$. This representation is unique up to an additive constant that is usually determined by imposing the condition that the function $g$ fixes the origin. The representation $f=h+\overline g$ is then unique and is called the \emph{canonical representation} of $f$.
\par
According to a theorem of Lewy \cite{Lewy}, a harmonic mapping $f=h+\overline g$ is locally univalent in $\D$ if and only if its \emph{Jacobian} $J_f=|h'|^2-|g'|^2$ is different from zero in the unit disk. Hence, every locally univalent harmonic mapping is either orientation preserving (if $J_f>0$ in $\D$) or orientation reversing (if $J_f<0$). Note that $f$ is orientation reversing if and only if $\overline f$ is orientation preserving. This trivial observation allows us to restrict ourselves to those cases when the locally univalent harmonic mappings considered preserve the orientation, so that $|h'|^2-|g'|^2>0$. Hence, the analytic function $h$ in the canonical representation of $f=h+\overline g$ is locally univalent and the \emph{dilatation} $\omega=g'/h'$ is analytic in $\D$ and maps the unit disk to itself.
\par\smallskip
It is plain that the harmonic mapping $f=h+\overline g$ is analytic if and only if the function $g$ is constant.

\subsection{Pre-Schwarzian and Schwarzian derivatives of harmonic mappings}
The \emph{harmonic pre-Schwarzian derivative} $P_H(f)$ and the \emph{harmonic Schwarzian derivative} $S_H(f)$ of an orientation preserving harmonic mapping $f=h+\overline g$ in the unit disk with dilatation $\omega=g'/h'$ were introduced in \cite{HM-Schwarzian}. These operators are defined, respectively, by the formulas
\[
P_H(f)=P(h)-\frac{\overline{\omega}\, \omega'}{1-|\omega|^2}
\]
and
\[
S_H(f)=S(h)+\frac{\overline \omega}{1-|\omega|^2} \left(\frac{h''}{h'}\, \omega'-\omega''\right)-\frac 32 \left(\frac{\overline \omega\, \omega'}{1-|\omega|^2}\right)^2\,,
\]
where $P(h)$ and $S(h)$ are the classical pre-Schwarzian and Schwarz\-ian derivatives of $h$.
\par
It is clear that when $f$ is analytic (so that its dilatation is constant), the harmonic pre-Schwarzian and Schwarzian derivatives of $f$ coincide with the classical definitions of the corresponding operators.
\par
The \emph{harmonic pre-Schwarzian} and \emph{Schwarzian norms} of the function $f$ are defined, respectively, by $\|P_H(f)\|=\sup_{z\in\D} |P_H(f)(z)|(1-|z|^2)$ and $\|S_H(f)\|=\sup_{z\in\D} |S_H(f)(z)|(1-|z|^2)^2$.
\par\smallskip
For further properties of the harmonic pre-Schwarzian and Schwarzian derivative operators and the motivation for their definition, see \cite{HM-Schwarzian}.

\par\smallskip
The Schwarzian operators $P_H$ and $S_H$ have proved to be useful to generalize classical results in the setting of analytic functions to the more general setting of harmonic mappings. This paper is another sample of their usefulness, as will become apparent in the proofs of our main results, Theorems~\ref{thm-main} and \ref{thm-main2} below.
\par
At this point, we mention explicitly the following criterion of univalence that generalizes the Nehari criterion stated above as well as the criterion for quasiconformal extension of locally univalent analytic functions due to Ahlfors and Weill \cite{AW}. The sharp value of the constant $\delta_0$ has still to be determined~\cite{HM-Nehari}.
\begin{other}\label{thm-nehariharmonic}
Let $f=h+\overline g$ be an orientation preserving harmonic mapping in $\D$. Then, there exists a positive real number $\delta_0$ such that if for all $z\in\D$
\[
\|S_H(f)\|=\sup_{z\in\D} |S_H(f)(z)|\, (1-|z|^2)^2\leq \delta_0\,,
\]
then $f$ is one-to-one in $\D$. Moreover, if $\|S_H(f)\|\leq \delta_0 t$ for some $t<1$, then $f$ has a quasiconformal extension to $\mathbb{C}\cup\{\infty\}$.
\end{other}
\par
The corresponding result, now in terms of the pre-Schwarzian derivative, is as follows (see \cite[Thm. 8]{HM-Schwarzian}). In this case, an extra-term involving the dilatation of the function $f$ must be taken into account. This extra-term is identically zero if $f$ is analytic, so that the next theorem is the generalization to the classical criterion of univalence due to Becker, Theorem~\ref{thm-BP}, to the cases when the functions considered are harmonic.
\begin{other}\label{thm-beckerharmonic}
Let $f=h+\overline g$ be an orientation preserving harmonic mapping in $\D$ with dilatation $\omega$. If for all $z\in\D$
\begin{equation}\label{eq-beckerharm}
|P_H(f)(z)|\, (1-|z|^2)+\frac{|\omega'(z)|\, (1-|z|^2)}{1-|\omega(z)|^2}\leq 1\,,
\end{equation}
then $f$ is univalent. The constant $1$ is sharp.
\end{other}
\par\smallskip
Criteria for quasiconformal extension of harmonic mappings in terms of the harmonic pre-Schwarzian derivative that extend the corresponding criteria in the analytic setting due to Becker \cite{Becker} and Ahlfors \cite{Ahlfors} can be found in \cite{HM-qc}.
\par\smallskip
We finish this section by pointing out the following remark that will be important later in this paper.
\par\smallskip
It is not difficult to check that if $f$ is an orientation preserving harmonic mapping and $\phi$ is an analytic function such that the composition $F=f\circ\phi$ is well defined, then $F$ is an orientation preserving harmonic mapping with dilatation $\omega_F=\omega\circ\phi$. Moreover, for all $z$ in the unit disk,
\begin{equation}\label{eq-chainrulepreS}
P_H(F)(z)=P_H(f) (\phi(z))\cdot \phi'(z)+\frac{\phi''(z)}{\phi'(z)}
\end{equation}
and
\begin{equation}\label{eq-chainruleSch}
S_H(F)(z)=S_H(f) (\phi(z))\cdot \left(\phi'(z)\right)^2+S(\phi)(z)\,.
\end{equation}
\subsection{Hyperbolic derivatives} Let $\omega$ be an analytic self-map of the unit disk (that is, $\omega$ is analytic in $\D$ and $|\/\omega(z)|<1$ for all $|z|<1$). The \emph{hyperbolic derivative} $\omega^*$ of such function $\omega$ is
\[
\omega^*(z)=\frac{\omega'(z)\, (1-|z|^2)}{1-|\omega(z)|^2}\,.
\]
\par
Notice that the second term in (\ref{eq-beckerharm}) coincides with $|\/\omega^*(z)|$.
\par\smallskip
The Schwarz-Pick lemma proves that $|\/\omega^*(z)|\leq 1$ for all $z$ in $\D$ and that equality holds at some point $z_0$ in the unit disk if and only if $\omega$ is an automorphism of $\D$. In this case, $|\/\omega^*|\equiv 1$.
\par
It is also easy to check that the \emph{chain rule} for the hyperbolic derivative holds: If $\omega$ and $\phi$ are two analytic self-maps of $\D$ and the composition $\omega\circ\phi$ is well-defined, then
\[
(\omega\circ\phi)^*(z)=\omega^*(\phi(z))\cdot \phi^*(z)\,.
\]
In particular,
\begin{equation}\label{eq-chainrulehypder}
|(\omega\circ\phi)^*(z)|\leq |\/\omega^*(\phi(z))|\,.
\end{equation}
\subsection{Valence of harmonic mappings} The zeros of a locally univalent harmonic mapping $f$ are isolated \cite[p. 8]{Duren-harm}. Just as in the analytic case, the \emph{valence} of such a harmonic function $f$ is defined by $\sup_{w\in\mathbb{C}} n(f,w)$, where $n(f,w)$ is the number of points $z\in\D$ (counting multiplicities) for which $f(z)=w$. The function $f$ is said to have \emph{bounded valence} if there exists a positive integer $N$ such that $\sup_{w\in\mathbb{C}} n(f,w)\leq N$.

\section{A criterion for bounded valence of harmonic mappings in terms of the pre-Schwarzian derivative}
We now state one of the two main theorems in this paper. It generalizes Theorem~\ref{thm-BP} to those cases when the function considered is just harmonic.
\begin{thm}\label{thm-main}
Let $f=h+\overline g$ be an orientation preserving harmonic mapping in the unit disk with dilatation $\omega$. If
\begin{equation}\label{eq-mainthm}
\limsup_{|z|\to 1} \left(|P_H(f)(z)|\, (1-|z|^2)+\frac{|\omega'(z)|\, (1-|z|^2)}{1-|\omega(z)|^2}\right)< 1\,,
\end{equation}
then there exists a positive integer $N$ such that $f$ takes every value at most $N$ times in $\D$.
\end{thm}
\par
It is possible  to show that if (\ref{eq-mainthm}) holds then all the analytic functions $\varphi_\lambda=h+\lambda g$, where $|\lambda|=1$, have bounded valence in the unit disk. However, we have not been able to prove directly that under the assumption that $\varphi_\lambda$ has bounded valence for all such $\lambda$, then $f$ has bounded valence too.
\par
The proof of our main theorem will follow similar arguments to those employed in the proof of Theorem~\ref{thm-BP}. However, the criterion of univalence needed in the case when the function $f$ is harmonic will be the one provided in Theorem~\ref{thm-beckerharmonic} instead of the classical criterion of univalence due to Becker. The following lemma will be needed to prove Theorem~\ref{thm-main}. We refer the reader to \cite[Lemmas 2.2 and 3.3]{Becker-Pom-2} (see also \cite{GGPPR}) for the details of the proof.
\begin{lem}\label{lem-aux}
Let $\rho\in (1/2, 1)$ and $\alpha>0$. Then, there exist a univalent analytic self-map $\psi$ of the unit disk and a positive integer $M$ such that
\[
\bigcup_{k=1}^M\ \left\{e^{\frac{2k\pi i}{M}} \psi(z)\ \colon\  z\in\D\right\}=\{\zeta\ \colon\ 2\rho-1<|\zeta|<1\}
\]
and
\[
\sup_{z\in\D} \left|\frac{\psi''(z)}{\psi'(z)}\right|\, (1-|z|^2)<\alpha\,.
\]
\end{lem}
\par\smallskip
We now prove Theorem~\ref{thm-main}.
\par\smallskip
\begin{pf}
By (\ref{eq-mainthm}), there exist  a real number $\rho$ with $1/2<\rho<1$ and $\beta<1$ such that
\begin{equation}\label{eq-mainpf1}
|P_H(f)(z)|\, (1-|z|^2)+\frac{|\omega'(z)|\, (1-|z|^2)}{1-|\omega(z)|^2}< \beta\,,\quad 2\rho-1<|z|<1\,.
\end{equation}
\par
Since the function $f$ is locally univalent and $|z|\leq 2\rho-1$ is compact, the function $f$ takes every value at most $L$ times in $|z|\leq 2\rho-1$.
\par
Let now $\psi$ be the univalent analytic self-map of the unit disk of Lemma~\ref{lem-aux} with $\alpha=(1-\beta)/2 >0$, so that for all positive integer $k\leq M$, the functions $\psi_k=e^{2k\pi i /M}\psi$ satisfy
\begin{equation}\label{eq-mainpf2}
\sup_{z\in\D} \left|\frac{\psi_k''(z)}{\psi_k'(z)}\right|\, (1-|z|^2)<\frac{1-\beta}{2}\,.
\end{equation}
\par
For any such value of $k$, define the functions $F_k=f\circ\psi_k$. These are orientation preserving harmonic mappings in the unit disk with dilatations $\omega_k=\omega\circ\psi_k$. Moreover, using (\ref{eq-chainrulepreS}), the triangle inequality, the Schwarz-Pick lemma, and  (\ref{eq-chainrulehypder}) we have that for all $|z|<1$,
\begin{align*}
\nonumber |P_H(F_k)(z)|\,& (1-|z|^2)+\frac{|\omega_k'(z)|\, (1-|z|^2)}{1-|\omega_k(z)|^2}\\
\nonumber &\leq |P_H(f)(\psi_k(z))|\, (1-|\psi_k(z)|^2)+ \left|\frac{\psi_k''(z)}{\psi_k'(z)}\right|\, (1-|z|^2)\\
\nonumber &+ \frac{|\omega'(\psi_k(z))|\, (1-|\psi_k(z)|^2)}{1-|\omega(\psi_k(z))|^2}\,.
\end{align*}
\par
Bearing in mind the fact that for all $z\in\D$ and all $k$ as above the modulus $|\psi_k(z)|>2\rho-1$, (\ref{eq-mainpf1}), and (\ref{eq-mainpf2}), we conclude
\[
|P_H(F_k)(z)|\, (1-|z|^2)+\frac{|\omega_k'(z)|\, (1-|z|^2)}{1-|\omega_k(z)|^2}\leq \beta+\frac{1-\beta}{2}= \frac{1+\beta}{2}<1\,.
\]
Hence, by Theorem~\ref{thm-beckerharmonic}, these functions $F_k=f\circ\psi_k$ are univalent in the unit disk. Since, by Lemma~\ref{lem-aux},
\[
\bigcup_{k=1}^M\ \left\{\psi_k(z)\colon z\in\D\right\}=\{\zeta\colon 2\rho-1<|\zeta|<1\}\,,
\]
it follows that $f$ takes every value at most $M$ times in $2\rho-1<|z|<1$, and we obtain that $f$ takes every value at most $N=L+M$ times in $\D$. This completes the proof.
\end{pf}
\section{Schwarzian derivative criterion for finite valence of harmonic mappings}

A direct consequence of the following lemma is that the Schwarzian derivative $S(\psi)$ defined by (\ref{eq-Schw-defn}) of the function $\psi$ from Lemma~\ref{lem-aux} will satisfy
\begin{equation}\label{eq-Schwarzianbound}
\sup_{z\in\D} |S(\psi)(z)|(1-|z|^2)^2 < 4\alpha+\frac{\alpha^2}{2}.
\end{equation}
Though the result is folklore (see, for instance, \cite[Proof of Lemma 10]{HKR}), we include the proof for the sake of completeness.
\begin{lem}\label{lem-Schwarzianbound}
Let $\psi$ be a locally univalent analytic function in the unit disk. Assume that
\[
\sup_{z\in\D} \left|\frac{\psi''(z)}{\psi'(z)}\right|\, (1-|z|^2)<\alpha\,.
\]
Then,
\[
\sup_{z\in\D} \left|\left(\frac{\psi''(z)}{\psi'(z)}\right)'\right|\, (1-|z|^2)^2<4\alpha\,.
\]
\end{lem}
\begin{pf}
In order to make the exposition more clear, let us use $\Psi$ to denote the analytic function $P(\psi)=\psi''/\psi'$.
\par
Given a fixed but arbitrary point $z\in\D$, let $r$ be the positive real number that satisfies $2r^2=1+|z|^2$. Hence,
\[
1-r^2=r^2-|z|^2=\frac{1-|z|^2}{2}\,.
\]

By hypotheses, for all $|\zeta|<1$,
\[
|\Psi(\zeta)|=\left|\frac{\psi''(\zeta)}{\psi'(\zeta)}\right|<\frac{\alpha}{1-|\zeta|^2}\,.
\]
The Cauchy and Poisson integral formulas now give
\begin{eqnarray*}
|\Psi'(z)|&=& \left|\frac{1}{2\pi i}\int_{|\zeta|=r} \frac{\Psi(\zeta)}{(\zeta-z)^2}\, d\zeta\right|\\
&<& \frac{\alpha}{1-r^2} \, \frac{1}{r^2-|z|^2} \frac{1}{2\pi}\int_0^{2\pi} \frac{r^2-|z|^2}{|re^{i\theta}-z|^2}\, d\theta\\
&=& \frac{\alpha}{1-r^2} \, \frac{1}{r^2-|z|^2}=\frac{4\alpha}{(1-|z|^2)^2}\,,
\end{eqnarray*}
which completes the proof.
\end{pf}
\par\smallskip
A criterion of bounded valence for harmonic mappings in the unit disk in terms or the harmonic Schwarzian derivative that generalizes Theorem~\ref{thm-GP} is as follows. The constant $\delta_0$ is equal to the one in Theorem~\ref{thm-nehariharmonic}.

\begin{thm}\label{thm-main2}
Let $f=h+\overline g$ be an orientation preserving harmonic mapping in the unit disk with dilatation $\omega$. If
\begin{equation}\label{eq-mainthm2}
\limsup_{|z|\to 1} |S_H(f)(z)|\, (1-|z|^2)^2< \delta_0\,,
\end{equation}
then $f$ has bounded valence in the unit disk.
\end{thm}
\begin{pf}
The argument of the proof is analogous to the one used to prove Theorem~\ref{thm-main}.
\par
Condition (\ref{eq-mainthm2}) implies that there exist a real number $\rho$ with $1/2<\rho<1$ and $\varepsilon>0$ such that
\begin{equation}\label{eq-mainpf3}
|S_H(f)(z)|\, (1-|z|^2)^2< \delta_0-\varepsilon\,,\quad 2\rho-1<|z|<1\,.
\end{equation}
\par
The function $f$ is locally univalent and $|z|\leq 2\rho-1$ is compact. Therefore~$f$ takes every value at most $L$ times in $|z|\leq 2\rho-1$.
\par
Consider the analytic self-map of the unit disk $\psi$ of Lemma~\ref{lem-aux} with $\alpha=\sqrt{16+2\varepsilon}-4$. Then, by Lemma~\ref{lem-Schwarzianbound}, we have that (\ref{eq-Schwarzianbound}) holds. Thus, for all positive integer $k\leq M$, the functions $\psi_k=e^{2k\pi i /M}\psi$ satisfy
\begin{equation}\label{eq-mainpf4}
\sup_{z\in\D} \left|S(\psi_k)(z)\right|\, (1-|z|^2)^2<4\alpha+\frac{\alpha^2}{2}=\varepsilon\,.
\end{equation}
\par
Using (\ref{eq-chainruleSch}), the triangle inequality, the Schwarz-Pick lemma, the fact that for all $z\in\D$ and all $k$ the modulus $|\psi_k(z)|>2\rho-1$, (\ref{eq-mainpf3}), and (\ref{eq-mainpf4}), we have that the functions $F_k=f\circ\psi_k$, $k=1, 2, \ldots, M$, will satisfy that for all $|z|<1$,
\begin{eqnarray}
\nonumber |S_H(F_k)(z)|\, 1-|z|^2)^2 &\leq& |S_H(f)(\psi_k(z))|\, (1-|\psi_k(z)|^2)^2\\
&+& \nonumber |S(\psi_k)(z)|(1-|z|^2)^2\\
&<& \nonumber \delta_0-\varepsilon + \varepsilon=\delta_0\,.
\end{eqnarray}
\par
Hence, by Theorem~\ref{thm-beckerharmonic}, these functions $F_k=f\circ\psi_k$ are univalent in the unit disk and, as in the proof of Theorem~\ref{thm-main}, it follows that $f$ takes every value at most $M$ times in $2\rho-1<|z|<1$. We then again obtain that $f$ takes every value at most $N=L+M$ times in $\D$.
\end{pf}




\end{document}